\documentclass[reqno]{amsart}
\usepackage{amsmath,amssymb}
\usepackage{hyperref}
\usepackage{graphicx}
\usepackage{fullpage}

\begin{document}
\title [Approximation of Functions of Lipschitz Class]{Error Bounds of Conjugate of a Periodic Signal by Almost Generalized N$\bf{\ddot{O}}$rlund Means }
\begin{center}
\author{Vishnu Narayan Mishra$^{1, 2}$,Vaishali Sonavane$^{1}$}
\end{center}
\address{Vishnu Narayan Mishra$^{1, 2}$ and Vaishali Sonavane$^{1}$
\newline $^{1}$ Applied Mathematics and Humanities Department,
Sardar Vallabhbhai National Institute of
Technology, Ichchhanath Mahadev Dumas Road, Surat, Surat-395 007 (Gujarat), India\\ 
 \newline $^{2}$ L1627 Awadh Puri Colony Beniganj, Phase – III, Opposite – Industrial Training Institute (I.T.I.), Ayodhya Main Road, Faizabad – 224 001 (Uttar Pradesh), India \newline}
 \email{vishnunarayanmishra@gmail.com; vishnu\_narayanmishra@yahoo.co.in and
vaishalisnvn@gmail.com}
\maketitle{\textbf{Abstract}}
In this paper, we determine the error bounds of conjugate signals between input periodic signals and processed output signals, whenever signals belong to $Lip\, (\alpha ,\, r)$ -class and as a processor we have taken almost generalized N$\ddot{o}$rlund  means using head bounded variation sequences and rest bounded variation sequences. The results obtained in this paper further extend several known results on linear operators. 
\noindent
\maketitle 
\section{\textbf{Introduction}}
\maketitle
\noindent A number of researchers Al-Saqabi et al.  \cite{2}, Liu and Srivastava  \cite{9}, Alzer et al.  \cite{10}, Mohapatra (\cite{22}-\cite{23}, Bor et al.  \cite{11}, Choi and Srivastava  \cite{15} have established interesting results in sequences and series using different linear summability operators. The degree of approximation of a function belonging to various classes using different linear summability operators has been determined by several investigators like Bernstein  \cite{22}, Alexits (\cite{7}-\cite{8}, Sahney and Goel  \cite{4}, Khan  \cite{13}, Mishra et al.(\cite{25}-\cite{26}. Chen and Hong  \cite{17} used Cesàro sum regularization technique for hyper singularity of dual integral equation. Summability of Fourier series is useful for engineering analysis. Sahney and Rao  \cite{3} and Khan  \cite{12}have studied the degree of approximation of functions belonging to $Lip\, (\alpha ,\, r)$ by $\left(N,p_{n} \right)$\& $\left(N,p,q\right)$ means respectively. Summability techniques were also applied on some engineering problems: for example, Chen and Jeng  \cite{16} implemented the $Ces\grave{a}ro$ sum of order (C, 1) and (C, 2), in order to accelerate the convergence rate to deal with the Gibbs phenomenon, for the dynamic response of a finite elastic body subjected to boundary traction. Recently, Mursaleen and Mohiuddine  \cite{21} discussed convergence methods for double sequences and their applications in various fields. Quereshi (\cite{18}-\cite{19} discussed the degree of approximation of function belonging to $Lip\, \alpha $ \&$Lip\left(\alpha ,\, r\right)$ by $\left(N,p_{n} \right)$ means of conjugate series of a Fourier series. Summation-Integral operators play an important role modeling various physical and biological processes. Analysis of signals or time functions is of great importance, because it conveys information or attributes of some phenomenon. The theory of linear summability operators has had an important influence on the development of mathematical systems theory. On the other hand, mathematical systems theory serves as a direct source of motivation and new techniques for the theory of linear summability operators and its applications. The engineers and scientists use properties of Fourier approximation for designing digital filters. Singh  \cite{23}, Pati  \cite{24} and Dikshit  \cite{14} have obtained the error bounds of conjugate signals by different summability methods. The purpose of this paper is to determine the error bounds of conjugate signals between input periodic signals and processed output signals, whenever signals belong to $Lip\, (\alpha ,\, r)$-class and as a processor by almost generalized N$\ddot{o}$rlund  means using under conditions that $\left(p_{n} \right)\, \in \, HBVS$and $\left(q_{n} \right)\, \in \, RBVS$. \\
\noindent Let \textit{$\sum _{n\, =\, 0}^{\infty }a_{n}  $}be a given infinite series with the sequence of partial sums$\left\{s_{n} \right\}$. Let p denotes the sequence $\{ p_{n} \} ,\, p_{-\, 1} \, =\, 0.$ For two sequences \textit{p} and \textit{q},
\[\begin{array}{l} {P_{n} :=p_{0} +p_{1} +\cdot \cdot {\rm \; }\cdot +p_{n} ,(P_{-1} =p_{-1} ={\rm \; }0),} \\ {Q_{n} :=q_{0} +q_{1} +\cdot \cdot {\rm \; }\cdot +q_{n} ,(Q_{-1} =q_{-1} ={\rm \; }0),} \end{array}\] 
the convolution  $(p*q)_{n} $  is defined by
\[R_{n} =(p*q)_{n} =\sum _{m\, =\, 0}^{n}p_{m} q_{n\, -\, m}  =\sum _{m\, =\, 0}^{n}p_{n\, -\, m} q_{m}  .\] 
It is obvious that
\[P_{n} :=(p*1)_{n} =\sum _{m\, =\, 0}^{n}p_{m}  \, \, {\rm and}\, \, Q_{n} :=(1*q)_{n} =\sum _{m\, =\, 0}^{n}q_{m}  \, =\, \sum _{m\, =\, 0}^{n}q_{n\, -\, m}  \, .\] 
When $(p*q)_{n} \, \ne \, 0$ for all \textit{$n$}, the generalized N$\ddot{o}$rlund  transform of the sequence $\left\{s_{n} \right\}$is the sequence $\left\{t_{\, n}^{\, p.q} \right\}$\textit{ }obtained by putting
\[t_{\, n}^{\, p.q} =\frac{1}{(p*q)_{n} } \sum _{m\, =\, 0}^{n}p_{n\, -\, m} \, q_{m}  \, s_{m} \, .\] 
If $t_{\, n}^{\, p.q} \to s\, \, \, {\rm as}\, \, \, n\to \infty ,$ then the sequence  $\left\{s_{n} \right\}$ is said to be summable to $s$ by generalized N$\ddot{o}$rlund  method $\left(N,p,q\right)$ and is denoted by $s_{n} \to s\left(N,p,q\right)$  \cite{5}.

\noindent The necessary and sufficient conditions for a $\left(N,p,q\right)$method to be regular are 
\[\sum _{m\, =\, 0}^{n}\left|p_{n\, -\, m} q_{m} \right| =O\, \left(\left|\, (p*q)_{n} \, \right|\right)\] 
 and 

\noindent $p_{n\, -\, m} =o\left(\left|(p*q)_{n} \right|\right),\, $as $\, n\to \infty ,$ for every fixed $m\ge 0$ for which $q_{m} \ne 0.$ 

\noindent The $\left(N,p,q\right)$method reduces to N$\ddot{o}$rlund  method $\left(N,p_{n} \right)$  if $q_{n} =1$  for all $n\, .$  The $\left(N,p,q\right)$method reduces to Riesz method$\left(\overline{N},q_{n} \right)$if  $p_{n} =1$  for all $n$. 

\noindent In the special case when $p_{n} =\left(\begin{array}{l} {n+\alpha -1} \\ {\, \, \, \, \alpha -1} \end{array}\right),\alpha >0,$ the method $\left(N,p_{n} \right)$reduces to the well-known method of summability $\left(C,\alpha \right).$ 

\noindent The particular case $p_{n} =\left(n+1\right)^{-1} $ of the N$\ddot{o}$rlund  mean is known as harmonic mean and is written as $\left(N,1/n+1\right).$ 

\noindent Let $f$be a$2\, \pi $- periodic signal (function) and Lebesgue integrable. The Fourier series of $f$ is given by   
\begin{equation} 
f\sim \frac{a_{0} }{2} +\sum _{n\, =\, 1}^{\infty }(a_{n} \cos \, n\, x+b_{n} \, sin\, n\, x) \, \equiv \, \sum _{n\, =\, 0}^{\infty }A_{n} \left(x\right)  
\end{equation} 
with $n^{th} $partial sum$s_{n} (f;\, x)$ called trigonometric polynomial of degree (or order) $n,$of the first (n+1) terms of the Fourier series of$f$.    

\noindent The conjugate series of Fourier series (1) is given by       
\begin{equation}
 \sum _{n\, =\, 1}^{\infty }(b_{n} \, \cos n\, x-a_{n} \sin n\, x) \equiv \, \sum _{n\, \, =\, \, 1}^{\infty }B_{n} \left(x\right) .                                                    \end{equation}
\noindent The conjugate function $\widetilde{f}(x)$ is defined for almost every $x$ by (see  \cite{6},definition 1.10).\\
\begin{align*}
\widetilde{f}(x)=-\frac{1}{2\pi } \int _{0}^{\pi }\psi \left(t\right)\, \cot {t\mathord{\left/ {\vphantom {t 2}} \right. \kern-\nulldelimiterspace} 2}  \, dt\, =\, \, \, \lim _{h\, \to \, 0} \left(-\frac{1}{2\pi } \int _{h}^{\pi }\psi \left(t\right)\, \cot {t\mathord{\left/ {\vphantom {t 2}} \right. \kern-\nulldelimiterspace} 2}  \, dt\, \right) 
\end{align*}
\noindent A signal (function)$f\in Lip\, \alpha $  if             
\begin{align*}
f(x+t)-f(x)={\rm O} \left(\left|\, t^{\, \alpha } \right|\right)\, \, \, \, \, {\rm for}\, \, \, 0<\alpha \le 1,\, \, t>0 
\end{align*}

\noindent and$f\in Lip\, (\alpha ,\, r),\, \, {\rm for}\, \, 0\le x\le 2\pi ,$ if
\begin{align*}
\left(\int _{0}^{2\, \pi }\left|\, f(x+t)-f(x)\, \right|\, ^{r} dx \right)^{{1\, \mathord{\left/ {\vphantom {1\,  \, r}} \right. \kern-\nulldelimiterspace} \, r} } ={\rm O} \left(\left|\, t\, \right|^{\, \alpha } \right),\, \, \, \, \, \, \, \, \, 0<\alpha \le 1,\, \, r\ge 1,\, \, t>0\, .
\end{align*} 
If we take $r\, \to \, \infty $then Lip($\alpha $, r)$=$ Lip $\alpha $.\\ 
\noindent $L_{\infty } $- norm of a function $f:\, R\, \to \, R$ is defined by $\left\| \, f\, \right\| _{\, \infty } \, = \sup \left\{ \left|\, f\left(x\right) \right|\, : x \in  R\right\}.$

\noindent $L_{r} $- norm of a function is defined by $\left\| f\right\| _{\, \, r} =\left(\int _{0}^{\, 2\, \, \pi }\left|f(x)\right|^{\, r} dx \right)^{{1\mathord{\left/ {\vphantom {1 r}} \right. \kern-\nulldelimiterspace} r} } ,\, \, \, 1\, \le \, \, r\, <\, \infty $.

\noindent The degree of approximation of a function$f:\, R\, \to \, R$ by trigonometric polynomial $t_{n} $of order $n$ under sup norm $\left\| \, \, \right\| _{\, \infty } $ is defined by ( \cite{1}) $r\, \to \, \infty $   
\begin{align*}
\left\| \, t_{n} \, -\, f\, \right\| _{\infty } \, =\, \sup \, \left\{\left|\, t_{n} \, (x)-f\left(x\right)\, \right|\, \, \, :\, \, x\, \in \, R\right\}  
\end{align*}            
and $E_{n} (f)$ of a function$f\in \, L_{\, r} $ is given by 
\begin{align*}
[E_{n} (f)=\mathop{\min }\limits_{n} \, \left\| \, \, t_{n} -f\, \right\| _{\, r} . 
\end{align*}
We say that the conjugate series is said to be almost generalized N$\ddot{o}$rlund  summable to the finite number $s,$ if 
\[\widetilde{t}_{n,\, r}^{p,\, q} =\frac{1}{R_{n} } \sum _{m\, \, =\, \, 0}^{n}p_{m} q_{n\, -\, m\, }  \tilde{s}_{m,\, r} \to s\, \, \, {\rm as}\, \, n\to \infty \, \] 
 uniformly with respect to $r,$where
\[\tilde{s}_{m,r} =\frac{1}{m+1} \sum _{j\, =\, r}^{m\, +\, r}\widetilde{s_{j} } .\] 
We note that $t_{n,\, \, r}^{p,\, \, q} $and $\widetilde{t}_{n,\, \, r}^{p,\, \, q} $are also trigonometric polynomials of degree (or order) n.

\noindent Here we defines two classes of sequences (see  \cite{20}).

\noindent A sequence $c=\left\{c_{n} \right\}$ of non-negative numbers is called Head Bounded Variation, or briefly $c\in HBVS{\rm ,}$if 
\begin{equation} 
\sum _{k\, =\, 0}^{n-1}\left|c_{k} -c_{k+1} \right| \le K\left(c\right)c_{n} ,\forall n\, \in \, N, 
\end{equation} 
or only for $n\le N$ if the sequence \textit{c} has only finite non-zero terms and the last non-zero term is $c_{N} .$\\
\noindent A sequence $c=\left\{c_{n} \right\}$ of non-negative numbers tending to zero is called of Rest Bounded Variation, or briefly $c\in R{\rm BVS,}$if 
\begin{equation} 
\sum _{k\, =\, \, n}^{\infty }\left|c_{k} -c_{k+1} \right| \le K\left(c\right)c_{n} ,\forall n\, \in \, N, 
\end{equation} 
where$K\left(c\right)$ is a constant depending only on \textit{c}. 

\noindent For an example of a sequence c $\in $ HBVS, we take $c$ to be any monotone increasing sequence. 

\noindent Then
\[\sum _{k\, =\, \, 0}^{n\, -\, 1}\left|c_{k} -c_{k\, +\, 1} \right| =\sum _{k\, =\, 0}^{n\, -\, 1}\left|c_{k\, +\, 1} -c_{k} \right| =c_{n} -c_{0} \le K\left(c\right)c_{n} ,\] 
  with $K\left(c\right)=1.$\\
\noindent For an RBVS sequence we take $c$ to be any monotone decreasing sequence with limit 0. Let $m$ be any positive integer with $m>n.$\\ Then
\[\sum _{k\, =\, \, n}^{m}\left|c_{k} -c_{k\, +\, 1} \right| =\sum _{k\, =\, n}^{m}\left|c_{k\, } -c_{k+1} \right| =c_{n} -c_{m+1} ,\] 
 and
\[\sum _{k\, =\, \, n}^{\infty }\left|c_{k} -c_{k\, +\, 1} \right| =c_{n} \mathop{\lim }\limits_{m} c_{m+1} =c_{n} -0=K\left(c\right)c_{n} ,\] 
 with, again, $K\left(c\right)=1.$

\noindent We shall use the following notations throughout the paper: 
\[\psi _{x} (t)\, =\, \psi \, (t)\, =\, f\, \left(x\, +\, t\right)\, -\, f\, (x\, -\, t)\, , \] 
\[\mathop{G_{n} }\limits^{\sim } (t)=\sum _{p\, =\, r}^{r\, +\, m}\frac{\cos \, (p+\, 1\, /2)\, t\cos pt/2}{\sin \, (t\, /\, 2)}  \, ,\] 
and M denotes a constant which may be different at each of its occurrence.\\
\section{\textbf{known result}}
\noindent Very recently Mishra et al.  \cite{27} determined the degree of approximation of a signal $f\, \in Lip(\alpha ,r),\, \, (r\ge 1)-$class by almost Riesz summability means of its Fourier series. Krasniqi  \cite{28} established the following theorem to estimate the error between the input signal $f\left(t\right)$and  signal obtained by passing through the almost generalized N$\ddot{o}$rlund  mean $t_{n,\, r}^{p,\, q} (f\left(t\right)\, ;\, x).$
\subsection{\textbf{Theorem }}
Let $\left(p_{n} \right)\, \in \, HBVS$ and $\left(q_{n} \right)\, \in \, RBVS$. If $f:{\mathbb R}\to {\mathbb R}$ is a $2\, \pi -$ periodic function, Lebesgue integrable and belonging to $Lip(\alpha ,r),\, \, (r\ge 1)-$ class, then the degree of approximation of the function $f$ by almost generalized N$\ddot{o}$rlund  means of its Fourier series $t_{n,\, r}^{p,\, q} (f\left(t\right)\, ;\, x)$ is given by\textbf{}
\begin{equation}
\left\| f(t)-t_{n,\, r}^{p,\, q} (f\left(t\right)\, ;\, x))\, \right\| _{\, r} =O\left(R_{\, n}^{{1\, \mathord{\left/ {\vphantom {1\,  \, r\, -\, \alpha }} \right. \kern-\nulldelimiterspace} \, r\, -\, \alpha } } \right),\forall n,
\end{equation}                                  
and  $\psi \left(t\right)$  satisfies the following conditions
\begin{equation} 
\left[\int _{0}^{{\pi \mathord{\left/ {\vphantom {\pi  R_{n} }} \right. \kern-\nulldelimiterspace} R_{n} } }\left(\frac{t\left|\psi )\left(t\right)\right|^{r} }{t^{\alpha } } \right)dt \right]^{{\raise0.7ex\hbox{$ 1 $}\!\mathord{\left/ {\vphantom {1 r}} \right. \kern-\nulldelimiterspace}\!\lower0.7ex\hbox{$ r $}} } =O\left(R_{n}^{-1} \right), 
\end{equation} 
\begin{equation}  
\left[\int _{{\pi \mathord{\left/ {\vphantom {\pi  R_{\, n} }} \right. \kern-\nulldelimiterspace} R_{\, n} } }^{\pi }\left(\frac{t^{-\delta } \left|\psi )\left(t\right)\right|^{r} }{t^{\alpha } } \right)dt \right]^{{\raise0.7ex\hbox{$ 1 $}\!\mathord{\left/ {\vphantom {1 r}} \right. \kern-\nulldelimiterspace}\!\lower0.7ex\hbox{$ r $}} } =O\left(R_{n}^{\delta -1} \right), 
\end{equation} 
where $\delta $ is a finite quantity, generalized N$\ddot{o}$rlund  means are regular  and $r+s=rs$ such that $1\le r\le \infty .$ 
\section{\textbf{ Main Result}} 

\noindent The object of this paper is to generalize the above result under much more general assumptions. We will measure the error between the input signal $\widetilde{f}\left(t\right)$and the processed output signal $\widetilde{t}_{n,\, r}^{p,\, q} (\widetilde{f}\left(t\right)\, ;\, x)$by establishing the following theorem: 

\subsection{\textbf{ Theorem}}
Let $f:{\mathbb R}\to {\mathbb R}$ be$2\, \pi -$periodic, integrable in the sense of Lebesgue and belonging to $Lip(\alpha ,r),\, \, (r\ge 1)-$class, then the degree of approximation of function $f$ by almost generalized N$\ddot{o}$rlund  means of its conjugate series of its Fourier series i.e. $\widetilde{t}_{n,\, r}^{p,\, q} (\widetilde{f}\left(t\right)\, ;\, x)$ is given by 
\begin{equation}
\left\| \widetilde{f}(t)-\widetilde{t}_{n,\, r}^{p,\, q} (f\left(t\right)\, ;\, x))\, \right\| _{\, r} =O\left(R_{n}^{{1\, \mathord{\left/ {\vphantom {1\,  \, r\, -\, \alpha }} \right. \kern-\nulldelimiterspace} \, r\, -\, \alpha } } \right),\forall n, 
\end{equation} 
$\left\{p_{n} \right\}\in HBVS\, \, {\rm and}\, \, \, \left\{q_{n} \right\}\in RBVS$ and $\psi \left(t\right)$ satisfies the following conditions
\begin{equation} 
\left[\int _{0}^{{\pi \mathord{\left/ {\vphantom {\pi  \, R_{n} }} \right. \kern-\nulldelimiterspace} \, R_{n} } }\left(\frac{\left|\psi )\left(t\right)\right|^{\, r} }{t^{\alpha } } \right)dt \right]^{{\raise0.7ex\hbox{$ 1 $}\!\mathord{\left/ {\vphantom {1 r}} \right. \kern-\nulldelimiterspace}\!\lower0.7ex\hbox{$ r $}} } =O\left(R_{\, n}^{\, -1} \right), 
\end{equation} 
\begin{equation} 
\left[\int _{{\pi \mathord{\left/ {\vphantom {\pi  R_{n} }} \right. \kern-\nulldelimiterspace} R_{n} } }^{\pi }\left(\frac{t^{-\delta } \left|\psi )\left(t\right)\right|\, ^{r} }{t^{\alpha } } \right)dt \right]^{{\raise0.7ex\hbox{$ 1 $}\!\mathord{\left/ {\vphantom {1 r}} \right. \kern-\nulldelimiterspace}\!\lower0.7ex\hbox{$ r $}} } =O\left(R_{\, n}^{\, \delta } \right), 
\end{equation} 
\noindent where $\delta $ is a finite quantity, generalized N$\ddot{o}$rlund  means are regular  and $r^{-1} +s^{-1} =1$ such that $1\le r\le \infty .$

\noindent \textbf{Proof :}

\noindent  Let us write 
\[\tilde{s}_{n} (f\, ;\, x)=\sum _{k\, \, =\, \, 1}^{n}B_{k} \left(x\right) \, \, \, {\rm and}\, \, \widetilde{t_{n} }\left(x\right)=\frac{1}{P_{n} } \sum _{k\, =\, \, 0}^{n}p_{n\, -\, k}  \tilde{s}_{k} \left(x\right).\] 
Then we have 
\begin{equation*}
{\tilde{s}_{k} (f\, ;\, x)-\tilde{f}(x)=\frac{1}{\pi } \int _{0}^{\pi }\psi _{x} (t)\frac{\cos \, (k\, +\, 1\, /2)\, t}{\sin \, (t\, /\, 2)} \,  dt} \\ 
=\frac{1}{\pi } \int _{0}^{\pi }\psi _{x} (t)\frac{\cos \, (k\, +\, 1\, /2)\, t}{\sin \, (t\, /\, 2)} \,  dt\, +\eta _{k} , 
\end{equation*}
where, by the Riemann-Lebesgue theorem,
\[\eta _{k} =\frac{1}{\pi } \int _{\delta }^{\pi }\psi (t)\frac{\cos \, (k\, +\, 1\, /2)\, t}{\sin \, (t\, /\, 2)} \,  dt\to 0\, \, \, \, {\rm as}\, \, \, \, k\to \infty ,\] 
and now
\begin{align*}
{\tilde{s}_{k,\, r} (f\, \left(t\right);\, \, x)-\tilde{f}(t)=\frac{2}{\pi \left(k+1\right)} \sum _{p=r}^{r+k}\int _{0}^{\pi }\psi (t)\frac{\cos \, (p+\, 1\, /2)\, t\cos pt/2}{\sin \, (t\, /\, 2)} \,  dt\, \, +\, \eta _{\, k\, ,\, \, r}  ,} \\ {\, \, \, \, \, \, \, \, \, \, \, \, \, \, \, \, \, \, \, \, \, \, \, \, \, \, \, \, \, \, \, \, } 
\end{align*}  
where, for almost generalized N$\ddot{o}$rlund  means of $\tilde{s}_{k,\, \, r} (f\, \left(t\right);\, x),$ we have 
\begin{align}
  \widetilde t_{n,\,r}^{p,\,q}(\widetilde f\left( t \right)\,;\,x)\, - \,\widetilde f(t) &= \frac{1}{{{R_n}}}\sum\limits_{m\, = \,0}^n {{p_m}{q_{n\, - \,m}}} \left\{ {{{\tilde s}_{m,\,r}}(f\,\left( t \right);\,x) - \tilde f(t)} \right\} \hfill\nonumber\\
  \,\,\,\,\,\,\,\,\,\,\,\,\,\,\,\,\,\,\,\,\,\,\,\,\,\,\,\,\,\,\,\,\,\,\,\,\,\,\,\,\,\, &= \frac{2}{{\pi {R_n}}}\int_0^\pi  {\psi (t)\sum\limits_{m\, = \,0}^n {\frac{{{p_m}{q_{n\, - \,m}}}}{{\left( {m + 1} \right)}}} \sum\limits_{p\, = \,r}^{r + m} {\frac{{\cos \,(p + \,1\,/2)\,t\cos pt/2}}{{\sin \,(t\,/\,2)}}} \,} dt\, + {\xi _{\,n\,,\,\,r}} \hfill\nonumber \\
  \,\,\,\,\,\,\,\,\,\,\,\,\,\,\,\,\,\,\,\,\,\,\,\,\,\,\,\,\,\,\,\,\,\,\,\,\,\,\,\,\,\,\, &= \frac{2}{{\pi {R_n}}}\int_0^\pi  {\psi (t)\sum\limits_{m\, = \,0}^n {\frac{{{p_m}{q_{n\, - \,m}}}}{{\left( {m + 1} \right)}}} \,} {\widetilde G_n}\left( t \right)dt\, + \,\,{\xi _{\,n\,,\,\,r}} \hfill \nonumber\\
  \,\,\,\,\,\,\,\,\,\,\,\,\,\,\,\,\,\,\,\,\,\,\,\,\,\,\,\,\,\,\,\,\,\,\,\,\,\,\,\,\,\, &= \frac{2}{{\pi {R_n}}}\left[ {\int\limits_0^{{\pi  \mathord{\left/
 {\vphantom {\pi  {{R_n}}}} \right.
 \kern-\nulldelimiterspace} {{R_n}}}} { + \int\limits_{{\pi  \mathord{\left/
 {\vphantom {\pi  {{R_n}}}} \right.
 \kern-\nulldelimiterspace} {{R_n}}}}^\pi  {} } } \right]\int_0^\pi  {\psi (t)\sum\limits_{m\, = \,0}^n {\frac{{{p_m}{q_{n\, - \,m}}}}{{\left( {m + 1} \right)}}} \,} {\widetilde G_n}\left( t \right)dt\, + \,{\xi _{\,n\,,\,\,r}} \hfill\nonumber\\
&= {I_1} + {I_2}\,({\text{say}}),
\end{align}
\noindent Using  H$\ddot{o}$lder inequality, $f\left(t\right)\in Lip(\alpha ,s)\, \Rightarrow \psi \left(t\right)\in Lip(\alpha ,s)$ on $\left[0,\pi \right]\, ,$ condition (9), inequalities
\begin{equation} 
\left(\sin \, t\, /\, 2\right)^{\, -\, 1} \, \le \, \pi \, /\, t\, ,\, \, {\rm for}\, \, 0\, <\, t\, \le \, \pi \, ,\left|\, \cos \, nt\, \right|\le 1. 
\end{equation} 
\begin{align}
\left| {{I_1}} \right| &\leqslant \frac{2}{{\pi {R_n}}}{\left[ {\int\limits_0^{{\pi  \mathord{\left/
 {\vphantom {\pi  {{R_n}}}} \right.
 \kern-\nulldelimiterspace} {{R_n}}}} {\left( {\frac{{{{\left| {\psi )\left( t \right)} \right|}^{\,r}}}}{{{t^\alpha }}}} \right)dt} } \right]^{{\raise0.7ex\hbox{$1$} \!\mathord{\left/
 {\vphantom {1 r}}\right.\kern-\nulldelimiterspace}
\!\lower0.7ex\hbox{$r$}}}}\,\,\,\,{\left[ {\int\limits_0^{{\pi  \mathord{\left/
 {\vphantom {\pi  {{R_n}}}} \right.
 \kern-\nulldelimiterspace} {{R_n}}}} {{{\left( {\frac{1}{{{t^{\, - \,\alpha }}}}\left| {\sum\limits_{m\, = \,0}^n {\frac{{{p_m}\,{q_{m - n}}}}{{m + 1}}} } \right|{{\widetilde G}_n}\left( t \right)} \right)}^{\,s}}dt} } \right]^{{\raise0.7ex\hbox{$1$} \!\mathord{\left/
 {\vphantom {1 s}}\right.\kern-\nulldelimiterspace}
\!\lower0.7ex\hbox{$s$}}}} \hfill \nonumber\\
   &= O\left( {R_n^{ - 1}} \right)\left[ {\int\limits_0^{{\pi  \mathord{\left/
 {\vphantom {\pi  {{R_n}}}} \right.
 \kern-\nulldelimiterspace} {{R_n}}}} {{{\left( {\frac{1}{{{t^{ - \alpha }}}}\,\left| {\sum\limits_{m\, = \,0}^n {\frac{{{p_m}\,{q_{m - n}}}}{{m\, + \,1}}} \,(m\, + \,1)\,{t^{ - 1}}\,} \right|} \right)}^{\,s}}dt} } \right] \hfill\nonumber \\ 
  \,& = O\left( {R_n^{ - 1}} \right){\left[ {\int\limits_0^{{\pi  \mathord{\left/
 {\vphantom {\pi  {{R_n}}}} \right.
 \kern-\nulldelimiterspace} {{R_n}}}} {R_n^{s\,\,}{t^{\left( {\alpha \, - \,1\,} \right)\,\,s}}dt} } \right]^{{{\,1} \mathord{\left/
 {\vphantom {{\,1} {\,s}}} \right.
 \kern-\nulldelimiterspace} {\,s}}}} = O\left( 1 \right)\,{\left[ {\int\limits_0^{{\pi  \mathord{\left/
 {\vphantom {\pi  {{R_n}}}} \right.
 \kern-\nulldelimiterspace} {{R_n}}}} {{t^{\left( {\alpha \, - 1} \right)\,\,s}}dt} } \right]^{{{\,1} \mathord{\left/
 {\vphantom {{\,1} s}} \right.
 \kern-\nulldelimiterspace} s}}} \hfill\nonumber\\
  & = O\left( 1 \right).O\left( {\frac{1}{{R_n^{\alpha \, - \,1\, + \,\frac{1}{s}}\,}}} \right) = O\left( {\frac{1}{{R_n^{\alpha \, - \,\left( {1\, - \,\frac{1}{s}} \right)}}}} \right)\, \hfill \nonumber\\ 
&=O\left(\frac{1}{R_{n}^{\alpha \, -\, \frac{1}{r} } } \right).
\end{align}  
\[\left|I_{2} \right|\le \frac{2}{\pi R_{n} } \left[\int _{{\pi \mathord{\left/ {\vphantom {\pi  R_{n} }} \right. \kern-\nulldelimiterspace} R_{n} } }^{\pi }\left(\frac{t^{-\delta } \left|\psi )\left(t\right)\right|^{\, r} }{t^{\alpha } } \right)dt \right]^{{\raise0.7ex\hbox{$ 1 $}\!\mathord{\left/ {\vphantom {1 r}} \right. \kern-\nulldelimiterspace}\!\lower0.7ex\hbox{$ r $}} } \, \, \, \, \left[\int _{{\pi \mathord{\left/ {\vphantom {\pi  R_{n} }} \right. \kern-\nulldelimiterspace} R_{n} } }^{\pi }\left(t^{\alpha \, \, +\, \, \delta } \left|\sum _{m\, =\, 0}^{n}\frac{p_{m} \, q_{m-n} }{m+1}  \right|\widetilde{G}_{n} \left(t\right)\right)^{\, s} dt \right]^{\, {\raise0.7ex\hbox{$ 1 $}\!\mathord{\left/ {\vphantom {1 s}} \right. \kern-\nulldelimiterspace}\!\lower0.7ex\hbox{$ s $}} } .\] 
Now, using the fact that $f\left(t\right)\in Lip(\alpha ,s)\, \Rightarrow \psi \left(t\right)\in Lip(\alpha ,s)$   on $\left[0,\pi \right],$ conditions (10),(12)and $r^{-1} +s^{-1} =1,$we obtain
\begin{align*}
\left| {{I_2}} \right| &\leqslant O\left( {R_n^{\delta \, - \,1}} \right)\,\,\,\,{\left[ {\int\limits_{{\pi  \mathord{\left/
 {\vphantom {\pi  {{R_n}}}} \right.
 \kern-\nulldelimiterspace} {{R_n}}}}^\pi  {{{\left( {\frac{{{t^{\alpha \, + \,\delta }}}}{{\sin t/2}}\sum\limits_{m\, = \,0}^n {\frac{{{p_m}{q_{m - n}}}}{{m + 1}}} \sum\limits_{p\, = \,r}^{r + m} {\cos \,(p + \,1\,/2)\,t\cos (pt/2)} } \right)}^s}dt} } \right]^{{\raise0.7ex\hbox{$1$} \!\mathord{\left/
 {\vphantom {1 s}}\right.\kern-\nulldelimiterspace}
\!\lower0.7ex\hbox{$s$}}}} \hfill \nonumber\\
  \,\,\,\,\,\,\,\, &= O\left( {R_n^{\delta \, - \,1}} \right)\,\,\,{\left[ {\int\limits_{{\pi  \mathord{\left/
 {\vphantom {\pi  {{R_n}}}} \right.
 \kern-\nulldelimiterspace} {{R_n}}}}^\pi  {{{\left( {\frac{{{t^{\alpha \, + \,\delta }}}}{{\sin t/2}}\sum\limits_{m\, = \,0}^n {{p_m}{q_{m - n}}} } \right)}^s}dt} } \right]^{{\raise0.7ex\hbox{$1$} \!\mathord{\left/
 {\vphantom {1 s}}\right.\kern-\nulldelimiterspace}
\!\lower0.7ex\hbox{$s$}}}}.\, \hfill \nonumber\\ 
\end{align*} 
Also, since $p=\left\{p_{n} \right\}\in \, {\rm HBVS,}$ then by (3), we obtain 
\begin{align}
 {p_m} - {p_n} &\leqslant \left| {{p_m} - {p_n}} \right| \leqslant \sum\limits_{k\, = \,m}^{n\, - \,1} {\left| {{p_k} - {p_{k\, + \,1}}} \right|}  \hfill\nonumber\\ &\leqslant \sum\limits_{k\, = \,0}^{n\, - \,1} {\left| {{p_k} - {p_{k\, + \,1}}} \right|}  \hfill\nonumber \\
& \leqslant K\left( p \right){p_n}, \hfill \nonumber\\
&\Rightarrow \,{p_m} \leqslant (K\left( p \right) + 1){p_n},\,\,\forall \,m \in \left[ {0,\,n} \right]. 
\end{align}
Also, since $q=\left\{q_{n} \right\}\in {\rm RBVS}$, then by (4), we obtain 
\begin{align}
{q_{n\, - \,m}} \leqslant \sum\limits_{k\, = \,m}^\infty  {\left| {{q_{n\, - \,k}} - {q_{n - \,k\, - 1}}} \right|}  \hfill\nonumber \\
  \,\,\,\,\,\,\,\,\,\, \leqslant \sum\limits_{k\, = \,0}^\infty  {\left| {{q_{n\, - \,k}} - {q_{n\, - \,k\, - \,1}}} \right|}  \hfill\nonumber\\ 
 \leqslant K\left( q \right){q_n},\forall \,m \in \left[ {0,\,n} \right].
\end{align} 
Now using (14)and (15), we have  
\begin{align}
\sum\limits_{m\, = \,0}^n {\left| {{p_m}{q_{m - n}}} \right|}  &\leqslant \sum\limits_{m\,\, = \,\,0}^n {(K\left( p \right) + 1){p_n}} \,K(q)\,{q_{\,n}} \hfill\nonumber\\
&= (K\left( p \right) + 1)K\left( q \right)\left( {n + 1} \right){p_n}{q_n} \hfill \nonumber\\ 
\,\, &= O\left( {{R_n}} \right).
\end{align}
Subsequently, we get
\begin{align} 
  \left| {{I_2}} \right|\, &= O\left( {R_n^\delta } \right)\,\,\,{\left[ {\int\limits_{{\pi  \mathord{\left/
 {\vphantom {\pi  {{R_n}}}} \right.
 \kern-\nulldelimiterspace} {{R_n}}}}^\pi  {{{\left( {{t^{\alpha \, + \,\delta \, - \,1}}} \right)}^s}dt} } \right]^{{\raise0.7ex\hbox{$1$} \!\mathord{\left/
 {\vphantom {1 s}}\right.\kern-\nulldelimiterspace}
\!\lower0.7ex\hbox{$s$}}}} \hfill\nonumber\\
  \,\,\,\,\,\,\,\, &= O\left( {R_n^\delta } \right)\,\,\,{\left[ {\left( {\frac{{{t^{(\alpha \, + \,\delta \, - \,1)\,s\, + \,1}}}}{{^{(\alpha \, + \,\delta \, - \,1)s\, + \,1}}}} \right)_{{\pi  \mathord{\left/
 {\vphantom {\pi  {{R_n}}}} \right.
 \kern-\nulldelimiterspace} {{R_n}}}}^\pi } \right]^{{\raise0.7ex\hbox{$1$} \!\mathord{\left/
 {\vphantom {1 s}}\right.\kern-\nulldelimiterspace}
\!\lower0.7ex\hbox{$s$}}}} \hfill\nonumber\\
  \,\,\,\,\,\,\,\,\nonumber\\ &= O\left( {R_n^\delta } \right)O\left( {R_n^{ - \,\alpha \, - \,\delta \, + \,1\, - \,1/s}} \right) \hfill\nonumber\\
  \,\,\,\,\,\,\,\, &= O\left( {R_n^{ - \,\alpha \, + \,1/r}} \right)\,\,\, \hfill\nonumber \\ 
&\, =O\left(\frac{1}{\left(R_{n}^{\alpha \, -\, 1/r} \right)} \right). 
\end{align}
 Using  (13)and (17) in (11) , we have
\[\left|\widetilde{f}(t)-\widetilde{t}_{n,\, r}^{p,\, q} (f\left(t\right)\, ;\, x))\right|=O\left(\frac{1}{R_{n}^{\alpha \, -\, 1/r} } \right).\] 
 Now, using $L_{r} $ -- norm, we obtain,
 \begin{align*}
{\left\| {\widetilde f(t) - \widetilde t_{n,\,r}^{p,\,q}(f\left( t \right)\,;\,x))\,} \right\|_r} &= {\left[ {\int\limits_0^{2\pi } {{{\left| {\widetilde f(t) - \widetilde t_{n,\,r}^{p,\,q}(f\left( t \right)\,;\,x))} \right|}^{\,r}}} dt} \right]^{{{\,1} \mathord{\left/
 {\vphantom {{\,1} r}} \right.
 \kern-\nulldelimiterspace} r}}} \hfill \\
  \,\,\,\,\,\,\,\,\,\,\,\,\,\,\,\,\,\,\,\,\,\,\,\,\,\,\,\,\,\,\,\,\,\,\,\,\,\,\,\,\,\,\,\,\,\,\,\,\, &= {\left[ {\int\limits_0^{2\pi } {O{{\left( {\frac{1}{{R_n^{\alpha \, - \,1/r}}}} \right)}^{\,r}}dt\,} } \right]^{{{\,1} \mathord{\left/
 {\vphantom {{\,1} r}} \right.
 \kern-\nulldelimiterspace} r}}}= O\left( {\frac{1}{{R_n^{\alpha \, - \,1/r}}}} \right). \hfill \\ 
 \end{align*}
This completes the proof of theorem 3.1.
\noindent 
\section{\textbf{Corollaries}} 
\noindent The result of our theorem is more general rather than the results of any other previous proved theorems, which will enrich the literate of summability theory of infinite series. The following corollaries can be derived from our main theorem:\textbf{}
\subsection{\textbf{Corollary}}
If  $f:{\mathbb R}\to {\mathbb R}$ be $2\, \pi -$ periodic, integrable in the sense of Lebesgue and belonging to $Lip\, \alpha $ class, then the degree of approximation of function $f$ by almost generalized N$\ddot{o}$rlund means of its conjugate series is given by     
\noindent 
\begin{equation}  
\left|\widetilde{f}(t)-\widetilde{t}_{n,\, r}^{p,\, q} (f\left(t\right)\, ;\, x))\, \right|_{\, \infty } =O\left(R_{n}^{-\alpha } \right),\forall n, 
\end{equation} 
${\rm where\; }\left\{p_{n} \right\}\in HBVS\, \, {\rm and}\, \, \left\{q_{n} \right\}\in RBVS\, \, \, $ and $\psi \left(t\right)$  satisfies (6), (7) and $r^{-1} +s^{-1} =1$ such that $1\le r\le \infty .$\\
\textbf{Proof} : Putting $r\to \infty $ in theorem 3.1 ,we have 
\noindent 
\[\left|\widetilde{f}(t)-\widetilde{t}_{n,\, r}^{p,\, q} (f\left(t\right)\, ;\, x))\right|_{\, \infty } =O\left(R_{n}^{-\alpha } \right).\] 
Finally, we find that
\[\left|\widetilde{f}(t)-\widetilde{t}_{n,\, r}^{p,\, q} (f\left(t\right)\, ;\, x))\right|\le \left|\widetilde{f}(t)-\widetilde{t}_{n,\, r}^{p,\, q} (f\left(t\right)\, ;\, x))\right|_{\, \infty } =O\left(R_{n}^{-\, \alpha } \right).\] 
For $\forall n\ge 0,$ we put $q_{n} =1\, \, \, {\rm or}\, \, p_{n} =1\, \, $ in corollary 4.1 , we obtain two cases\\
\noindent  \textbf{1.} If $f:{\mathbb R}\to {\mathbb R}$ be $2\, \pi -$ periodic, integrable in the sense of Lebesgue and belonging to $Lip\, \alpha $ class, then the degree of approximation of function $f$ by almost Riesz means of its conjugate series                          
\begin{equation} 
\left|\widetilde{f}(t)-\widetilde{t}_{n,\, r}^{p,\, q} (f\left(t\right)\, ;\, x))\, \right|=O\left(P_{n}^{-\, \alpha } \right),\forall n, 
\end{equation} 
${\rm where\; }\left\{p_{n} \right\}\in HBVS\, $ and $\psi \left(t\right)$  satisfies the following conditions
\begin{equation}  
\left[\int _{0}^{{\pi \mathord{\left/ {\vphantom {\pi  P_{n} }} \right. \kern-\nulldelimiterspace} P_{n} } }\left(\frac{\left|\psi )\left(t\right)\right|^{\, r} }{t^{\alpha } } \right)dt \right]^{\, {\raise0.7ex\hbox{$ 1 $}\!\mathord{\left/ {\vphantom {1 r}} \right. \kern-\nulldelimiterspace}\!\lower0.7ex\hbox{$ r $}} } =O\left(P_{n}^{-1} \right), 
\end{equation}                    
\begin{equation}
\left[\int _{{\pi \mathord{\left/ {\vphantom {\pi  P_{n} }} \right. \kern-\nulldelimiterspace} P_{n} } }^{\pi }\left(\frac{t^{-\, \delta } \left|\psi )\left(t\right)\right|^{\, r} }{t^{\, \alpha } } \right)dt \right]^{\, {\raise0.7ex\hbox{$ 1 $}\!\mathord{\left/ {\vphantom {1 r}} \right. \kern-\nulldelimiterspace}\!\lower0.7ex\hbox{$ r $}} } =O\left(P_{n}^{\, \delta } \right), 
\end{equation} 
 where $\delta $ is a finite quantity, Riesz means are regular  and $r^{-1} +s^{-1} =1$ such that $1\le r\le \infty .$\\
\noindent  \textbf{2.} If $f:{\mathbb R}\to {\mathbb R}$ be $2\, \pi -$periodic, integrable in the sense of Lebesgue and belonging to $Lip\, \alpha $ class, then the degree of approximation of function $f$ by almost N$\ddot{o}$rlund  means of its conjugate series   
\noindent 
\begin{equation} 
\left|\widetilde{f}(t)-\widetilde{t}_{n,\, r}^{p,\, q} (f\left(t\right)\, ;\, x))\, \right|=O\left(Q_{n}^{-\, \alpha } \right),\forall n,\, \,  
\end{equation} 
   $\left\{q_{n} \right\}\in RBVS\, \, \, $and $\psi \left(t\right)$satisfies the following conditions
\begin{equation} 
\left[\int _{0}^{{\pi \mathord{\left/ {\vphantom {\pi  Q_{n} }} \right. \kern-\nulldelimiterspace} Q_{n} } }\left(\frac{\left|\psi )\left(t\right)\right|^{\, r} }{t^{\, \alpha } } \right)dt \right]^{\, {\raise0.7ex\hbox{$ 1 $}\!\mathord{\left/ {\vphantom {1 r}} \right. \kern-\nulldelimiterspace}\!\lower0.7ex\hbox{$ r $}} } =O\left(Q_{n}^{-1} \right), 
\end{equation} 
\begin{equation}  
\left[\int _{{\pi \mathord{\left/ {\vphantom {\pi  Q_{n} }} \right. \kern-\nulldelimiterspace} Q_{n} } }^{\pi }\left(\frac{t^{-\delta } \left|\psi )\left(t\right)\right|^{\, r} }{t^{\, \alpha } } \right)dt \right]^{\, {\raise0.7ex\hbox{$ 1 $}\!\mathord{\left/ {\vphantom {1 r}} \right. \kern-\nulldelimiterspace}\!\lower0.7ex\hbox{$ r $}} } =O\left(Q_{n}^{\delta } \right), 
\end{equation} 
where $\delta $ is a finite quantity, means are regular  and $r^{-1} +s^{-1} =1$ such that $1\le r\le \infty .$
\noindent Let $q_{n} =1$ in theorem 3.1 $\forall n\ge 0,$  we get
\subsection{\textbf{Corollary}} If $f:{\mathbb R}\to {\mathbb R}$ be $2\, \pi -$periodic, integrable in the sense of Lebesgue and belonging to $Lip\left(\alpha \, ,r\right)\, ,\, (r\, \ge \, 1)$ class, then the degree of approximation of function $f$ by almost Riesz means of its conjugate series is given by    
\noindent 
\begin{equation} 
\left|\widetilde{f}(t)-\widetilde{t}_{n,\, r}^{p,\, q} (f\left(t\right)\, ;\, x))\, \right|=O\left(P_{n}^{{\, 1\, \mathord{\left/ {\vphantom {\, 1\,  \, r-\alpha }} \right. \kern-\nulldelimiterspace} \, r-\alpha } } \right),\forall n,\, \, \, \,  
\end{equation} 
${\rm where\; }\left\{p_{n} \right\}\in HBVS\, ,$$\psi \left(t\right)$ satisfies the following conditions (20) and (21), $\delta $ is a finite quantity, Riesz means are regular  and $r^{-1} +s^{-1} =1$ such that $1\le r\le \infty .$\\
\noindent Let $p_{n} =1$ in theorem 3.1 $\forall n\ge 0,$  we get
\subsection{\textbf{Corollary}} If $f:{\mathbb R}\to {\mathbb R}$ be $2\, \pi -$periodic, integrable in the sense of Lebesgue and belonging to $Lip\left(\alpha \, ,p\right)$class, then the degree of approximation of function $f$ by almost N$\ddot{o}$rlund  means of its conjugate series        
\begin{equation}  
\left|\widetilde{f}(t)-\widetilde{t}_{n,\, r}^{p,\, q} (f\left(t\right)\, ;\, x))\, \right|=O\left(Q_{n}^{{\, 1\mathord{\left/ {\vphantom {\, 1 \, r-\alpha }} \right. \kern-\nulldelimiterspace} \, r-\alpha } } \right),\forall n,\, \, \, \,  
\end{equation} 
${\rm where\; }\left\{q_{n} \right\}$ in RBVS, $\psi \left(t\right)$ satisfies the following conditions (23) and (24,$\delta $ is a finite quantity, N$\ddot{o}$rlund means are regular  and $r^{-1} +s^{-1} =1$ such that $1\le r\le \infty$ .


\begin{thebibliography}{30}
\bibitem{1} A. Zygmund, Trigonometric Series, V. Warszawa (1935), 40-41.
\bibitem{2} B.N. Al-Saqabi, S.L. Kalla and H.M. Srivastava, A Certain Family of Infinite Series Associated with Digamma Functions, J. Math. Anal. Appl. 159,(1991) 361-372 .
\bibitem{3} B.N. Sahney and V.V.G. Rao, Error bounds in the approximation of functions, Bull. Aust. Math. Soc., 6(1972), 11-18.
\bibitem{4} B.N. Sahney and D.S. Goel, On degree of approximation of continuous functions, RanchiUniv. Math. J., 4(1973), 50-53.
\bibitem{5} D. Borwein, On product of sequences, J. London Math. Soc., Vol. 33, (1958), 352 -- 357.
\bibitem{6} E. Hille and J. D. Tamarkin, On the summability of Fourier series. I, Trans.  of the Amer. Math. Soc. 34 (4) (1932), 757-783.
\bibitem{7} G. Alexits, Uber die Annanherung einer stetigin Function durch die Cesaroschen Mittel inhrer Fourier-riche, Mathe Annal, 100(1928), 264-277.
\bibitem{8} G. Alexits, Convergence Problems of Orthogonal Series, Pergamon Elmsford, New York, NY, USA, 1961.
\bibitem{9} G.D. Liu and H. M. Srivastava, Explicit formulas for the  N$\ddot{o}$rlund polynomials $B_{n}^{(x)} $ and $b_{n}^{(x)} ,$ Computers \& Mathematics with Applications 51 (2006), 1377-1384.
\bibitem{10} H. Alzer, D. Karayannakis and H.M. Srivastava, Series representations for some mathematical constants, J. Math. Anal. Appl. 320 (2006), 145--162.
\bibitem{11} H. Bor, H. M. Srivastava and W. T. Sulaiman, A new application of certain generalized power increasing sequences, Filomat 26:4 (2012), 871--879 doi 10.2298/FIL1204871B.
\bibitem{12} H.H. Khan, On the degree of approximation of a functions belonging to the class $Lip\left(\alpha ,p\right),$ Indian  J. Pure Appl. Math. 5 (1974), 132--136.
\bibitem{13} H.H. Khan, A note on a theorem of Izumi, Commun. Fac. Maths. Ankara 31 (1982), 123--127.
\bibitem{14} H.P. Dikshit, The Nörlund summability of the conjugate series of a Fourier series, Rendiconti del Circolo Matematico di Palermo 1 May--August 1962,~Volume 11,~Issue 2\underbar{,}~pp 217-224.
\bibitem{15} J.  Choi and H. M. Srivastava, Certain Classes of Series Involving the Zeta Function. J. Math. Anal. Appl.  231(1991), 91--117.
\bibitem{16} J.T. Chen and Y.S. Jeng, Dual series representation and its applications to a string subjected to support motions, Advances in Engineering Software, vol. 27, no. 3, (1996), pp. 227--238.
\bibitem{17} J.T. Chen and H.K. Hong, Review of dual boundary element methods with emphasis on hypersingular integrals and divergent series, Applied Mechanics Reviews, ASME, 52 (1999), 17-33.
\bibitem{18} K. Qureshi, On the degree of approximation of functions belonging to the Lipschitz Class by means of a conjugate series, India J. Pure Appl. Math., Vol. 12, No. 9, (1981), 1120-1123.
\bibitem{19} K. Qureshi, On the degree of approximation of functions belonging to the class Lip ($\alpha $, r), by means of a conjugate series, Indian J. Pure Appl., Vol. 13, No. (5), (1982), 560--563.
\bibitem{20} L. Leindler, On the uniform convergence and boundedness of a certain class of sine series, Anal. Math. 27, (2001), 279 - 285.
\bibitem{21} M. Mursaleen and S. A. Mohiuddine, Convergence Methods For Double Sequences and Applications, Springer, 2014.
\bibitem{22} R.N. Mohapatra, Quantitative results on almost convergence of a sequence of positive linear operators, J. Approx. Theory, Vol. 20, Issue 3, 1977, 239-250.
\bibitem{23} R.N. Mohapatra, Functions of class Lip ($\alpha $,p) and their Taylor mean, J. Approx. Theory, Vol. 45, Issue 4, 1985, 363-374.
\bibitem{24} S. Bernstein, Surl Order de la Meilleure approximation desfunctions continues pardes polynomes de degree, donne's Memories Acad. Roy. Belyique (2), 4(1912), 1-104.
\bibitem{25} T. Singh, N$\ddot{o}$rlund summability of Fourier series and its conjugate series, Annali di Matematica Pura ed Applicata, Volume 64,~Issue 1,~(1964), pp 123-132.
\bibitem{26} T. Pati, On the Absolute N$\ddot{o}$rlund Summability of the conjugate series of a Fourier Series\textbf{,} J. London Math. Soc.~(1963)~s1-38(1):~204-214.
\bibitem{27} V.N. Mishra, K. Khatri and L.N. Mishra, Approximation of functions belonging to Lip($\xi $(t), r) class by (N, p${}_{n}$) (E, q) summability of conjugate series of Fourier series, Journal of Inequalities and Applications, (2012) 2012:296. doi: 10.1186/1029-242X-2012-296.
\bibitem{28} V.N. Mishra , K. Khatri  and L.N. Mishra, Using linear operators to approximate signals of ${\rm Lip}\left(\alpha ,p\right),\left(p\ge {\rm 1}\right)$-class, Filomat 27:2 (2013), 353--363. doi 10.2298/FIL1302353M.
\bibitem{29} V.N. Mishra, H.H. Khan, I.A. Khan and L.N. Mishra, On the degree of approximation of signals $Lip\left(\alpha \, ,r\right),\, \left(r\ge 1\right)$ class by almost Riesz means of its Fourier series, Journal of Classical Analysis 4 (2014), 79-87.
\bibitem{30} X.Z. Krasniqi, Applications of some classes of sequences on approximation of functions (signals) by almost generalized N$\ddot{o}$rlund means of their Fourier series, accepted in Journal of Classical Analysis (2015), in press.   
\end{thebibliography}
\end{document}